\documentclass[reprint,aps,pra,superscriptaddress]{revtex4-2}
\usepackage{amsmath}
\usepackage{amsfonts}
\usepackage{amsthm}
\usepackage{amssymb}
\usepackage{algorithm}
\usepackage{algorithmicx}
\usepackage[algo2e]{algorithm2e}
\usepackage{algcompatible}
\usepackage{algpseudocode}
\usepackage{bm}
\usepackage{nicematrix}
\usepackage{graphicx}
\usepackage{subfigure}
\begin{document}
\title{A Split Fast Fourier Transform Algorithm for Block Toeplitz Matrix-Vector Multiplication}
\author{Alexandre Siron} 
\author{Sean Molesky}
\affiliation{Department of Engineering Physics, Polytechnique Montr\'{e}al, Montr\'{e}al, Qu\'{e}bec H3T 1J4, CAN}
\begin{abstract}
    \noindent
    Numeric modeling of electromagnetics and acoustics frequently entails matrix-vector multiplication with block Toeplitz structure. 
    When the corresponding block Toeplitz matrix is not highly sparse, e.g. when considering the electromagnetic Green function in a spatial basis, such calculations are often carried out by performing a multilevel embedding that gives the matrix a fully circulant form.  
    While this transformation allows the associated matrix-vector multiplication to be computed via Fast Fourier Transforms (FFTs) and diagonal multiplication, generally leading to dramatic performance improvements compared to naive multiplication, it also adds unnecessary information that increases memory consumption and reduces computational efficiency. 
    As an improvement, we propose a lazy embedding, eager projection, algorithm that for dimensionality $d$, asymptotically reduces the number of needed computations $\propto d/ \left(2 - 2^{-d+1}\right)$ and peak memory usage $\propto 2/\left((d+1)2^{-d} + 1\right)$, generally, and $\propto\left(2^{d} + 1\right)/\left(d +2\right)$ for a fully symmetric or skew-symmetric systems.
    The structure of the algorithm suggests several simple approaches for parallelization of large block Toeplitz matrix-vector products across multiple devices and adds flexibility in memory and task management.
\end{abstract}
\maketitle
\begin{figure*}[htp]
\centering
\includegraphics[width=1\linewidth]{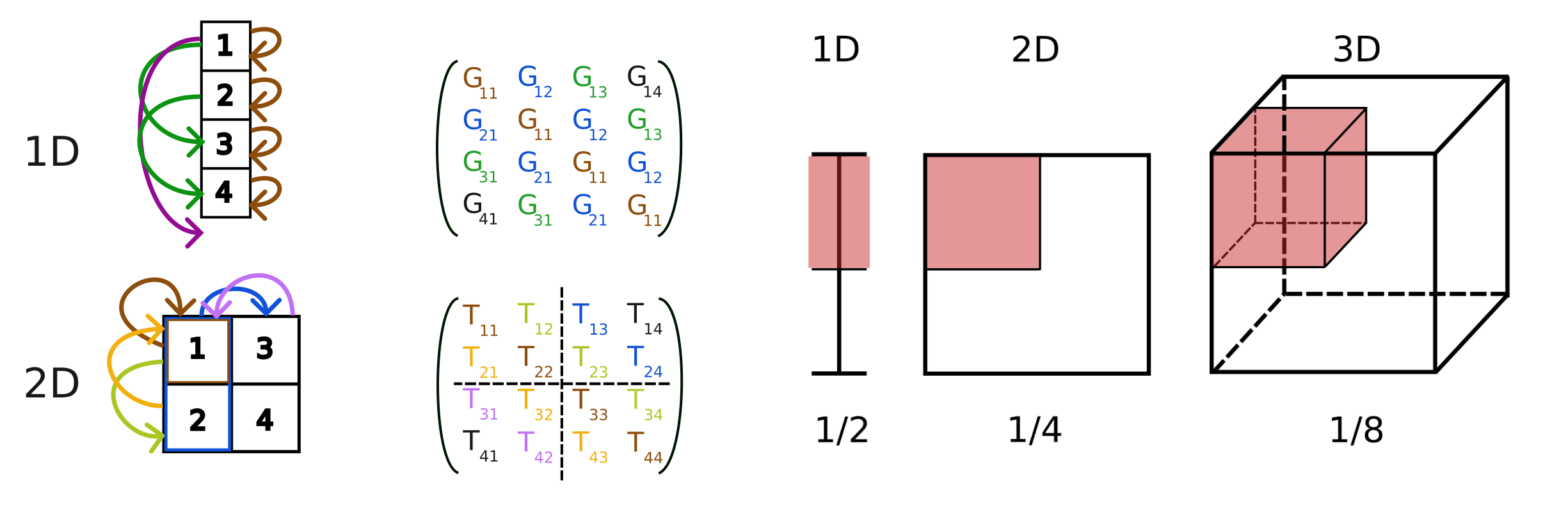}
\vspace{-36pts}
\caption{\label{fig:transinv}Left: 
\textbf{Translational invariance and block Toeplitz structure.} 
The schematic displays the distribution of interaction coefficients for translationally invariant physics of a pair of simple grids. 
0-step (self-interactions) are labeled in brown, vertical 1-step interactions in yellow and light green, horizontal 1-step interactions in blue and pink, 2-step interactions in darker green, and 3-step interactions in purple. 
Right: 
\textbf{Block Toeplitz embedding and data padding.} 
The cartoon illustrates the dilution of input vector data in the standard circulant embedding procedure as the dimensionality (level) of block Toeplitz structure increases. 
Input vector data is is indicated as the coloured portion of each sub-figure.}
\end{figure*}

Translationally invariant physics, when considered in a spatial basis, leads to matrix representations with (multi-level) block Toeplitz matrix structure, Fig.~\ref{fig:transinv}. 
Through scattering theory---wherein a complex physical system is described in terms of its free space properties and ``scattering'' effects caused by previously unaccounted for interactions \cite{ref15}---the need to efficiently compute block Toeplitz matrix-vector products accordingly arises in a variety of physical contexts~\cite{ref1}.
For example, block Toeplitz matrices appear frequently in simulation and device design problems for electromagnetic~\cite{ref2,ref3} and acoustic waves~\cite{ref4}, as well as the analysis of multi-dimensional discrete random processes~\cite{ref5}. 
More directly, in many of these settings (e.g. the simulation of nanophotonic devices utilizing Green functions \cite{ref16}) a sort of scattering phenomena is efficiently described by a highly sparse operation (e.g. the added interaction is local in space) and the primary bottle-neck of iterative methods is the speed of block Toeplitz matrix-vector multiplication~\cite{ref19}. 

The standard method of computing block Toeplitz matrix-vector products---motivated a number of perspectives~\cite{ref6,ref7}---is to perform a circulant embedding for each level of the Toeplitz structure. 
That is, the Toeplitz matrix-vector method of ``extending the diagonal bands'',
\begin{equation}
    \underbrace{
        \begin{pNiceArray}{cc}
            t_{11} & t_{12}\\
            t_{21} & t_{11}\\
        \end{pNiceArray}}_{n \times n~\text{Toeplitz}} \rightarrow
    \underbrace{    
        \begin{pNiceArray}{cc|cc}
          t_{11} & t_{12} & s_0 & t_{21} \\
          t_{21} & t_{11} & t_{12} & s_0 \\
          \hline
          s_0 & t_{21} & t_{11} & t_{12} \\
          t_{12} & s_0 & t_{21} & t_{11}
        \end{pNiceArray}}_{2n \times 2n~\text{circulant}},
        \nonumber
\end{equation}
is carried out iteratively.
Through this procedure, the full matrix acquires a circulant form, and multiplication with the  resulting embedded vector is carried out via the FFT approach to convolution \cite{ref7}. 
Notably, each embedding doubles the size of the corresponding Toeplitz level, and hence overall size of the system, by the inclusion of additional zero coefficients: taking $d$ as the number of levels, the embedding results in a vector space that is $2^d$ times as large.
These padding coefficients occupy an increasingly large fraction of the total system information as the dimension (number of levels) of the Toeplitz structure increases [Fig. \ref{fig:transinv}], leading to inefficient memory utilization and computation.

Two major categories of alternatives have been previously suggested to this procedure. 
The first class focuses on how block Toeplitz matrix-vector multiplication happens mechanistically, examining how successive coefficient reorderings may be used to minimize the number of arithmetic operations~\cite{ref8}. 
These methods, to the best of our understanding, have yet to offer substantial speed improvements over the usual FFT approach for large systems, and can not be simply implemented with optimized functions from standard libraries.
The second class focuses on the acceleration of block Toeplitz (or Toeplitz) matrix-vector mulplications through the use of matrix or tensor product decompositions, e.g. via Vandermonde matrices~\cite{ref9}, QTT tensor formatting~\cite{ref10}, or trigonometric transformation~\cite{ref11}. 
While extremely powerful in certain settings \cite{ref10}, these approach generally require some degree of approximation, are more computationally expensive than circulant embedding for arbitrary vectors, and rely on specialized code bases. 
Similarly, proposals to improve contemporary FFT algorithms, either under specific conditions~\cite{ref12} or by providing a functional equivalent~\cite{ref13}, are, at present, only beneficial for sparse vectors.  
However, if progress is made in this area, it should be directly applicable to the algorithm presented below. 

Here, recognizing that the data layout of a block Toeplitz structure allows for the ordering of embedding, Fourier transformation, and projection to be partially interchanged, we present an algorithm that performs block Toeplitz matrix-vector multiplication by dividing each dimension into a pair of branches. 
In this ``divide and conquer'' approach, circulant embeddings are replaced by the creation of phase modified vector copies containing even and odd Fourier coefficients (II.A). 
This change enables a lazy evaluation strategy that can substantially reduce memory pressure, $\propto\left(2^{d} + 1\right)/\left(d +2\right)$ in the case of electromagnetics, and operational complexity. 
Moreover, the partitioned nature of the algorithm presents several possibilities for tunable parallelization to increase computational speed, or handle large block Toeplitz systems. 
An accompanying implementation used for all presented results, written in Julia, is available from https://github.com/alsirc/SplitFFT\_lazyEmbed.
As the procedure is straightforward (II.B), and FFT dominated, similar performance should be simply achievable in any major scientific computing language. 

\begin{figure*}[htp]
\centering
\includegraphics[width=1\linewidth]{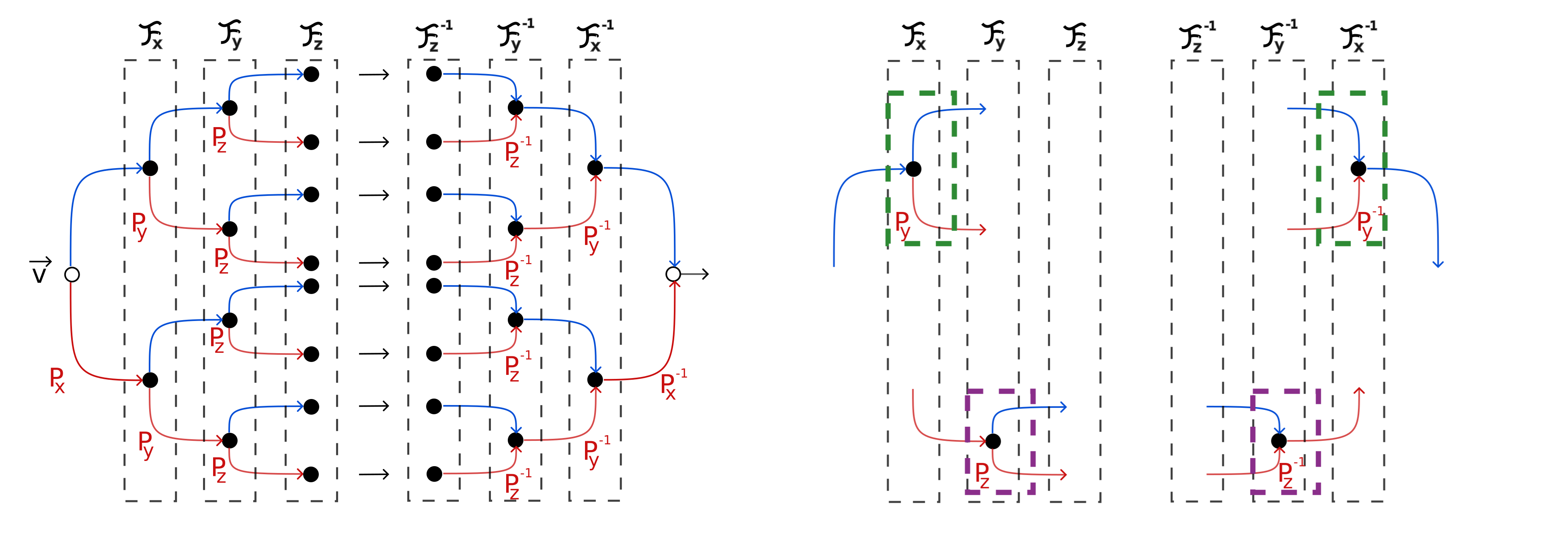}
\vspace{-24pts}
\caption{\label{fig:treebranch}Left: 
\textbf{Branching form of proposed algorithm.} 
The left panel depicts the branching tree of data flow in 3D to perform block Toeplitz matrix-vector multiplication in the proposed algorithm, Alg.~\ref{alg:cap}. 
Black dots are FFT tranformations, $P$ labelled arrows represent odd (phase modified) Fourier coefficients, and black horizontal arrows are diagonal multiplications with block Toeplitz matrix data. 
Right: 
\textbf{Branch execution in recursive implementation.} 
The schematic shows the particular parts of the data flow that are controlled by a given branch (bId) in reference to Alg.~\ref{alg:cap}.
Two different levels of \textit{toeMulBrn} are displayed, green and purple boxes.
Each level launches two function calls, and then merges the retruned results. 
Control is then returned to the calling level.}
\end{figure*}

\section{Algorithm}
Conceptually, the proposed algorithm is based on an observation about how the padding elements of a given level of embedding are incorporated into subsequent FFTs.
Working from the finest level of Toeplitz structure (numbers) to the coarsest (outer matrix blocks), the added zeros are only mixed with the input vector information once the FFT of the particular level has been completed. 
Analogously, as soon as the inverse FFT of a particular level has been carried out, half of the elements no longer influence the final matrix-vector product: they are obviated by subsequent projection, Fig.~\ref{fig:transinv} right. 
As such, by postponing each embedding step until it is necessary, and performing each projection as soon as possible, required memory and computation can be reduced. 

The properties described under subsection A indicate how lazy (postponed) embedding and eager (as soon as possible) projection are achieved by a division of the input vector at each level into ``original'' and phase-modified ``child'' copies that carry even and odd Fourier coefficients. 
Via this branching structure, the overall matrix-vector operation is partitioned into smaller calculations, which can be leveraged to reduce memory pressure---information is only calculated once it is needed to continue evaluation of the topmost (all-even) branch Fig.~\ref{fig:treebranch}. 
The structure also offers possibilities for additional parallelization. 
For instance, in a multi-device setting, the linearity of the Fourier transform could be used to perform parallel matrix-vector multiplication in a number of ways: a branch may either communicate and merge, decreasing overall computation, or continue to some other predefined termination point, decreasing communication overhead. 
To handle larger systems, the decomposition structure of the Fourier transform could be used to further divide the number of coefficients treated at any one time.
That is, $\mathbf{v}$ could be pre-partitioned into nearly independent even and odd parts.

\subsection{Splitting and merging of vector coefficients}
At each level, splitting into even and odd indices is accomplished using the following property.
Let $\mathbf{v}$ be a vector of size $n$. 
Define $\mathbf{\Tilde{v}}$ as $\mathbf{\Tilde{v}} = \left[\mathbf{v}, 0\right]$,
a vector of size $2n$, and 
$P$ as the diagonal phase shift operator, $P = \text{diag}\left[1,\ldots, e^{i \pi k/n}, \ldots, -1\right]$. 
$\mathcal{F}(\mathbf{v})$ and $\mathcal{F}(P\mathbf{v})$ are then, respectively, the even and odd coefficients of $\mathcal{F}(\mathbf{\Tilde{v}})$, where $\mathcal{F}$ denotes the Fourier transform. 
\begin{proof}
    Let $v_k$ denote coefficients of the vectors $\mathbf{v}$, and $\tilde{v}_{k}$ the coefficients of $\mathbf{\Tilde{v}}$.
    \begin{equation}
        \mathcal{F}(\mathbf{\Tilde{v}}) = \left( \frac{1}{2n} \sum_{k=0}^{2n-2}\Tilde{v}_k \Tilde{\omega}^{lk} \right)_{l=0,...n-1}\!\!\!\!,
        \nonumber
    \end{equation}
    where $\Tilde{\omega} = \exp(-\frac{i\pi}{n})$, and
    \begin{equation}
        \mathcal{F}(\mathbf{v}) = \left(\frac{1}{n} \sum_{k=0}^{n-1} v_k\omega^{lk}\right)_{l=0,...,n-1} \!\!\!\!\!\!\!\!\!\! = \left(\frac{1}{n} \sum_{k=0}^{n-1} v_k\Tilde{\omega^{2lk}}\right)_{l=0,...,n-1}\!\!\!\!,
        \nonumber
    \end{equation}
    where $\omega = \exp(-\frac{2i\pi}{n})$.
    Similarly,
    \begin{equation}
        \mathcal{F}(P\mathbf{v}) = \frac{1}{\sqrt{n}}\left(\sum_{k=0}^{n-1} \frac{\omega^{lk}}{\sqrt{n}}\exp\left(i\frac{l\pi}{n}\right)\right)_{l=0,...,n-1}\!\!\!\!\!\!\!\!\!\!\times \mathbf{v}
        \nonumber
    \end{equation}
    \begin{equation}
        = \frac{1}{n}(\sum_{k=o}^{n-1} \Tilde{\omega}^{(2l+1)k} v_k)_{l=0,...n-1}.
        \nonumber
    \end{equation}
    Therefore, 
    $
        \nonumber
        \mathcal{F}(\Tilde{\mathbf{v}}) = \mathcal{F}(\Tilde{\mathbf{v}}_{\text{e}}) + \mathcal{F}(\Tilde{\mathbf{v}}_{\text{o}}) = \mathcal{F}({\mathbf{v}}) + \mathcal{F}(P\mathbf{v}).
    $
\end{proof}
\noindent
Through this separation, the circulant embedding operation at a given level of the block Toeplitz structure is transformed into a branching operation, with each branch retaining the size of the initial vector.
\\ \\
After the diagonal multiplication step, branches are merged in reverse order---the projection counterpart to embedding. 
Using the equivalence presented above, accounting for the difference in relative normalization factors, merging is accomplished via the formula 
\begin{equation}
    \Tilde{\mathbf{v}} = \mathcal{F}^{-1}\left(\frac{1}{2}(\mathcal{F}(\mathbf{v}_{\text{even}}) + \overline{P} \mathcal{F}(\mathbf{v}_{\text{odd}}))\right).
    \nonumber
\end{equation} 
\subsection{Procedure}
Making use of the above properties, the proposed algorithm follows the pseudo-code presented under \emph{Algorithm 1}. 
Additional details on the Julia implementation available via https://github.com/alsirc/SplitFFT\_lazyEmbed are given in the appendix.
\begin{algorithm}[H]
    \caption{\!\!: block Toeplitz multiplication}\label{alg:cap}
    \textbf{Require} $d \geq 0$\\
    \textbf{Ensure} $y = Tv$\\
    $y \gets v$\\
    \If{$d > 0$}{\!\\
        \Comment{FFT along the $d^{th}$ dimension}\\
        $v \gets FFT_{d}(v)$\\
    }
    \eIf{$d < d_{\text{max}}$}{\!\\
        \Comment{create phase shifted vector (splitting)}\\
        $v_{\text{child}} \gets sptBrn(v)$ \\
        \Comment{launch next branches (recursive call)}\\
        $v \gets toeMulBrn(T, d+1, bId, v)$\\
        $v_{\text{child}} \gets toeMulBrn(T, d+1, nxtId(T, d), v_{\text{child}})$\\
        \Comment{merge $v$ (even coeffs.) with $v_{\text{child}}$ (odd coeffs.)}\\
        $v \gets mrgBrn(T, v, d+1, P, v_{\text{child}})$\\}{\!\\
    \Comment{diagonal multiplication of $v$ with Toeplitz data}\\
    $v \gets mulBrn(v, T[bId])$}
    \If{$d > 0$}{\!\\
        \Comment{inverse FFT along the $d^{th}$ dimension}\\
    \STATE $v \gets iFFT_d(v)$
    }
\end{algorithm}
\vspace{-12pts}
The protocol proceeds recursively by nested calls to the \textit{toeMulBrn} function, which organizes the splitting, merging, and possible diagonal multiplication operations that occur on a specific level, or dimension, $d$. 
As visualized in Fig.\ref{fig:treebranch}, the first task completed by \textit{toeMulBrn} is to compute the FFT along the current level. 
\textit{toeMulBrn} then creates a phase shifted copy of the transformed vector by applying the \textit{sptBrn} function, and launches two progeny branches, respectively containing even and odd Fourier coefficients. 
Overall coordination is ensured by passing each of these calls a branch identifier, \textit{bId}: the even branch inherits the identifier of its caller, while the odd branch is assigned a new identifier by \textit{nxtId}. 
The function then waits for these calls to return, merges their results, using \textit{mrgBrn}, and performs an inverse FFT, \textit{IFFT}, before returning.
If the maximal level of block Toeplitz structure is reached, $d = d_{\text{max}}$ then \textit{toeMulBrn} launches no further calls, and instead computes the diagonal multiplication of its given vector with the appropriate Toeplitz data, $T\left[bId\right]$, via \textit{mulBrn}. 
\textit{toeMulBrn} then computes the inverse FFT, $\textit{IFFT}$, and returns.
\subsection{Operational complexity and memory usage}
Consider a $d$-dimensional vector $v$ of total size $s = n^d$ where $d > 0$. 
As before, let $\Tilde{v}$ be the fully embedded image of $v$ (along all dimensions)---total length of $\tilde{v}$ is $2^d s$.
The complexity of a FFT as used in a level of Alg.~\ref{alg:cap} is $m_2 \log_2(m_1)$ where $m_1$ is the total length of the vector to be transformed and $m_2$ is the product of its sizes along the dimensions to be transformed~\cite{ref12}. 
The complexity of the standard circulant embedding method is therefore
\begin{equation}
    C_{\text{embed}} =  \underbrace{2^{d+1} s \log_2(2^d s)}_{\text{FFTs}} + \underbrace{2^d s}_{\text{multiplication}}.
    \nonumber
\end{equation}
At the $l$\textsuperscript{th} level, $1 < l < d$, of Alg.~\ref{alg:cap} there are $2^l$ forward FFTs and $2^{l}$ inverse FFTs, each along one dimension. 
Including the additional diagonal multiplications involved with the phase shifting and multiplication with the Toeplitz data, the operational complexity of Alg.~\ref{alg:cap} is thus
\begin{flalign}
    C_{\text{split}} &= \underbrace{2 \sum^{d}_{l=1} 2^l s\log_2(n)}_{\text{FFTs}} + \underbrace{2^d s}_{\text{multiplication}} + \underbrace{2 \sum^{d-1}_{l=0} 2^l s}_{\text{phase shift}} &
    \nonumber \\
    &= 2(2^d - 1)s(2\log_2(n) + 1) + 2^d s
    \nonumber
\end{flalign}
Take $R_c = C_{\text{embed}}/C_{\text{split}}$ to be operational complexity ratio between the two methods. 
Combining the above,
\begin{flalign}
    R_c &= \frac{2^{d+1}s\log_2(2^d s) + 2^d s}{2(2^d - 1)s(2\log_2(n) + 1) + 2^d s} &
    \nonumber \\
    &= \frac{d\log_2(2n) + 1}{(1 - 2^{-d})(2\log_2(n) + 1) + 1}, &
    \label{CompRatio}
\end{flalign}
so that $R_c\rightarrow d/ \left(2 - 2^{-d+1}\right)$ as $n \rightarrow \infty$.

In terms of memory usage, the standard circulant method requires an overall doubling in vector size for each level (dimension) of block Toeplitz structure. 
The memory required for both the Toeplitz data and input vector is, consequently, $M_{v,\text{embed}} = 2^d s = M_{T,\text{embed}}$, with $s$ as above. 
For Alg.~\ref{alg:cap}, a lazy evaluation strategy leads peak memory usage to occur when the first branch (top of Fig.~\ref{fig:treebranch}) is processed up to multiplication with its corresponding Toeplitz data, leading to a peak allocation of $M_{v,\text{split}} = (d+1)s$ elements. 
For the Toeplitz information, the full embedded vector of Fourier transformed data is split into $2^d$ parts of $s$ coefficients each, so $M_{T,\text{split}} = 2^d s$. 
Comparing these results, the peak memory ratio $R_m$ of the two approaches, in general, is
\begin{flalign}
    R_m &= \displaystyle \frac{M_{v,\text{embed}} + M_{T, \text{embed}}}{M_{v,\text{split}} + M_{T, \text{split}}} 
    =\frac{2^{d+1} s}{\left(d + 1\right)s + 2^d s} &
    \nonumber \\
    &= \frac{2}{\left(d + 1\right)2^{-d} + 1}. &
    \label{MemRatio}
\end{flalign}
Because electromagnetics, acoustics, and quatum mechanics typically exhibit reciprocity between sources and receivers \cite{ref14}, there are a number of settings of physical interest wherein system matrices are both block Toeplitz and block symmetric (resp. block skew-symmetric). 
In such cases, the Fourier transform of the vector containing the generating elements of the associated Toeplitz matrix is mirror symmetric (resp. skew-symmetric) about the middle index of each dimension, and thus can be fully stored in a vector of size $s$. 
In such cases, the peak memory ratio between the two methods becomes
\begin{equation}
    R_{m, sym} = \frac{2^{d} + 1}{\left(d + 2\right)}.
    \label{MemRatioSym}
\end{equation}
\noindent
More precisely, if $\mathbf{v}$ is a vector with coefficients $(v_k)_{k=0,n-1}$ such that $v_i = v_{n-i}$ for $i = 1,..n-1$, as happens in the case of a symmetric block Toplitz matrix, then its Fourier coefficents $\mathbf{f}$ exhibit the symmetry $f_i = f_{n/2 + i}$. 
Hence, only half the coefficients must be stored. 
Inductively, an equivalent reduction occurs at each level, and single vector's worth of Toeplitz data is required. 
The antisymmetric case is analogous as $v_i = -v_{n-i}$ implies that $f_{i+n/2} = -f_i - 2f_0$.\\ \\
\begin{figure*}[htp]
\centering
\subfigure[]{\includegraphics[width=0.45\linewidth]{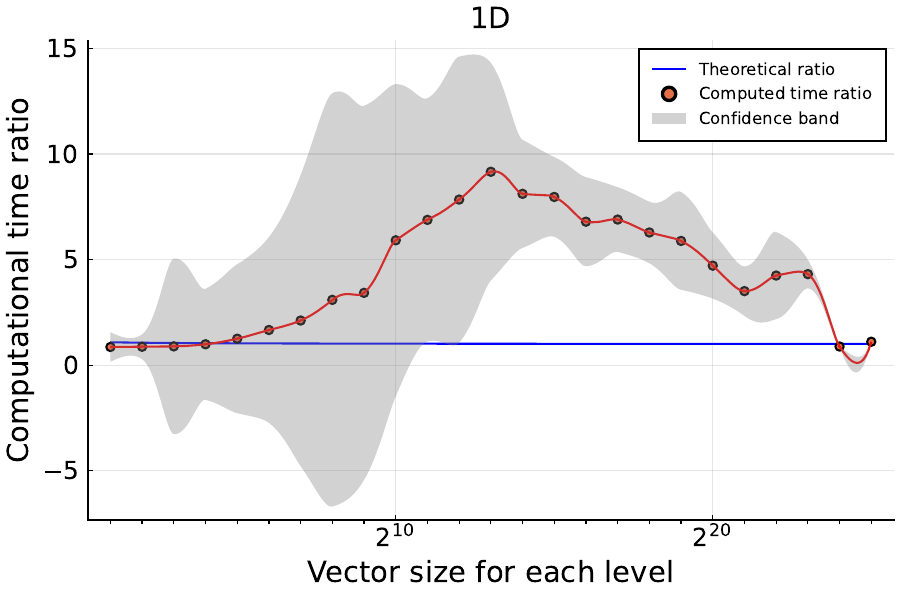}}
\subfigure[]{\includegraphics[width=0.45\linewidth]{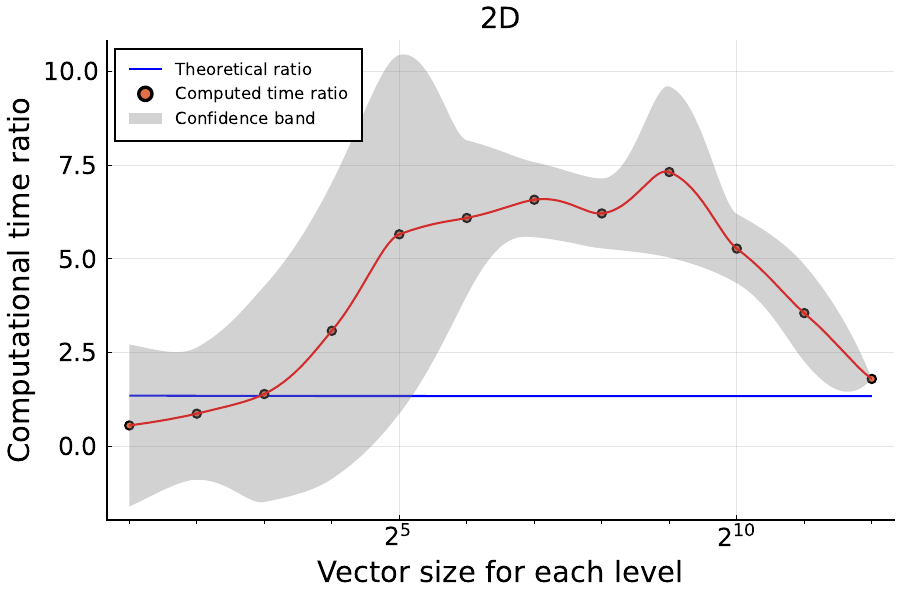}}
\subfigure[]{\includegraphics[width=0.45\linewidth]{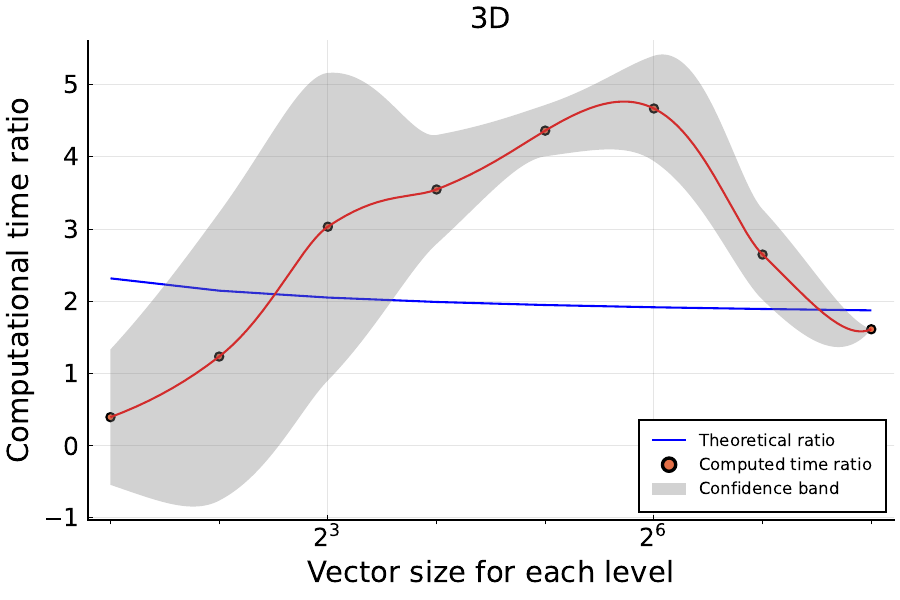}}
\subfigure[]{\includegraphics[width=0.45\linewidth]{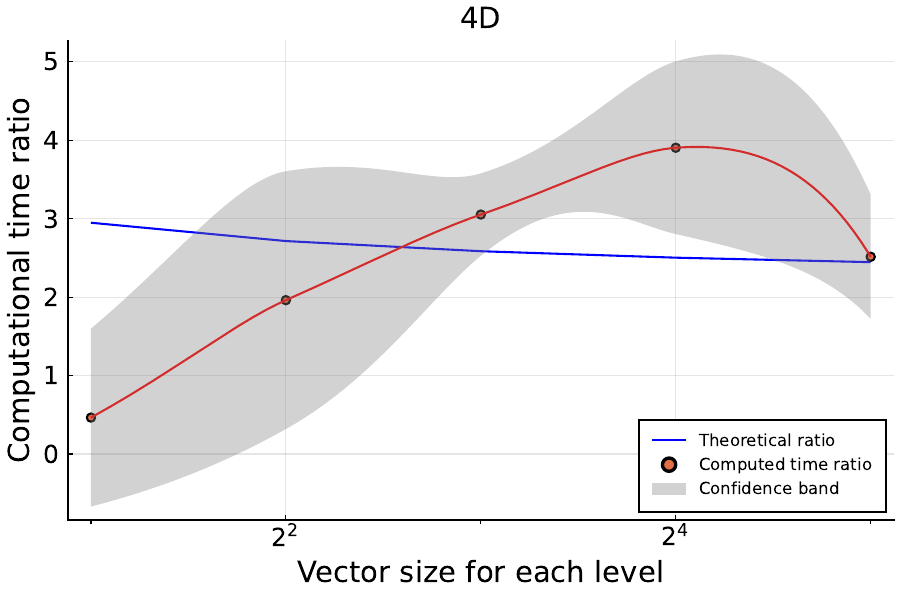}}
\caption{\textbf{Comparsion of wall clock and operational complexity ratios.} 
The figure depicts the wall clock time ratio of the Julia implementation of Alg.~\ref{alg:cap} (https://github.com/alsirc/SplitFFT\_lazyEmbed) compared to standard circulant embedding for vectors with lengths corresponding to powers of $2$.
In moving from panel (a) to (d), the dimensionality of block Toeplitz structure is increased from $1$ (Toeplitz matrix) to $4$.}
\label{fig:plot1-4d}
\end{figure*}

\textbf{Example of coefficient symmetry for two levels}\\
$$
    \begin{bmatrix}
        f^{00} & f^{01} & f^{02} & f^{03} & 0 & f^{03} & f^{02} & f^{01} \\
        f^{10} & f^{11} & f^{12} & f^{13} & 0 & f^{13} & f^{12} & f^{11} \\
        f^{20} & f^{21} & f^{22} & f^{23} & 0 & f^{23} & f^{22} & f^{21} \\
        f^{30} & f^{31} & f^{32} & f^{33} & 0 & f^{33} & f^{32} & f^{31} \\
        0 & 0 & 0 & 0 & 0 & 0 & 0 & 0 \\
        f^{30} & f^{31} & f^{32} & f^{33} & 0 & f^{33} & f^{32} & f^{31} \\
        f^{20} & f^{21} & f^{22} & f^{23} & 0 & f^{23} & f^{22} & f^{21} \\
        f^{10} & f^{11} & f^{12} & f^{13} & 0 & f^{13} & f^{12} & V^{11}
    \end{bmatrix}
$$
\begin{table}[h!]
    \centering
    \begin{tabular}{||c c c c c c||} 
 \hline
 Dimensions & 2 & 3 & 4 & 5 & 6 \\ [0.25ex] 
 \hline\hline
 $R_c$ & 1.33 & 1.71 & 2.13 & 2.58 & 3.05 \\ [0.65ex]
 \hline
 $R_m$ & 1.14 & 1.33 & 1.52 & 1.68 & 1.80 \\ [0.65ex]
 \hline
 $R{m,\text{sym}}$ & 1.25 & 1.80 & 2.83 & 4.71 & 8.13 \\ [0.65ex]
 \hline
\end{tabular}
    \caption{Asymptiotic theoretic complexity, Eq.~\ref{CompRatio}, and peak memory ratios, Eq.~\ref{MemRatio} and \ref{MemRatioSym}.}
    \label{tab:complexity}
\end{table}
\section{Results}
As analytic complexity calculations must necessarily overlook many, potentially influential, practical implementation details, we directly measure the relative performance of the split and standard circulant embedding algorithms by performing block Toeplitz matrix-vector multiplication for various dimension  and vector sizes. 

For convenience, the length of each dimension is assumed to be consistent, but this is not a limitation of either method (or code).
Our presented results focus on lengths corresponding to powers of $2$, as FFTs are most efficient for these lengths. 
Plots of cases where the sizes are products of primes and plots of consecutive natural numbers are available in the appendix. 
Versions of the presented plots with uncertainty bars corresponding to the statistical variance of the time computation are also in the appendix.
All results were obtained through local computation on a laptop with 16 GiB of RAM, and an AMD® Ryzen 7 5800hs processor.
Although the Julia implementation can be run on GPU, for simplicity, only parallelized, multi-threaded, CPU computations are compared.
The coefficients used for the matrix vectors and the input vector are complex numbers with random real and imaginary part with mean value equal to 0 and a parametrizable variance. 
The number of points for each plot is limited by the memory capacity of the machine used for the computation.

The results obtained do not clearly show the tendencies predicted by operational complexity calculations.
Namely, Alg.~\ref{alg:cap} is generally faster than expected, and only converges to the predictions of Eq.~\eqref{CompRatio} as the memory limit of the test machine is reached.
The broad uncertainty observed (see appendix) suggests that the statistical data captured by BenchmarkTools.jl is independent and uncorrelated.
These differences in measured performance, compared to Eq.~\eqref{CompRatio}, and large variance of results, may be caused by several factors. 
First, all FFTs are computed using the FFTW package \cite{ref17}, which utilizes different methods depending on the structure of the vector to be transformed, and has an internal thread management system that is not controlled by Julia. 
For a fixed input vector, the two methods perform FFTs of different sizes, and smaller sizes are known to be comparatively faster than expected based theoretical complexity~\cite{ref18}.
Second, the standard circulant embedding method employs full thread parallelization on fully embedded vectors at all steps. 
Conversely, Alg.~\ref{alg:cap} consistently applies thread parallelization at the size of the input vector and makes greater use of the Julia task scheduler. 
For small numbers of tasks, and small loops, launch overhead becomes a meaningful percentage of the benchmarks. 
Parallelization of small loops can also lead to highly variable memory writes and reads, which could explain the spread of measured ratios.

\section{Outlook}
In this article, we have described how the standard circulant embedding method for multi-level block Toeplitz matrix-vector multiplication can be improved by reordering embeddings and projections with respect to their associated Fourier transformation level. 
Concretely, by utilizing a lazy embedding / eager projection strategy, for $d$ levels of Toeplitz structure, we have shown that operational complexity can be reduced by a factor of $d/ \left(2 - 2^{-d+1}\right)$, and peak memory usage by a factor of $2/\left(\left(d+1\right)2^{-d} +1\right)$. 
As the size of the matrix-vector product tends to infinity, for $3$ levels of structure, these ratios tend to $12/7$ and $4/3$ respectively. 
For a fully symmetric (resp. skew-symmetric) matrix peak memory reduction improves to a factor of $\left(2^{d}+1\right)/\left(d+2\right)$---$9/5$ for $d=3$. 
In simulation, the wall clock time ratio between the two methods is better than theoretical expectations in almost all cases. 
The proposed algorithm also suggests several way in which additional parallelization could be introduced to handle the large scattering problems that occur in acoustics and electromagnetics. 
In fact, a slightly modified version of the (GPU enabled) Julia implementation is already being used in the GilaElectromagnetics package (https://github.com/moleskySean/GilaElectromagnetics.jl), currently under development by the authors. 
The algorithm relies only on (highly optimized) standard library functions, and should be easily adaptable to almost any relevant context.

\section{Acknowledgments}
The authors thank Justin Cardona and Paul Virally for their careful reading and feedback on this article, and Paul Virally for his assistance in configuring various packages. 
This work was supported by the Québec Ministry of Economy of Innovation through the Québec Quantum Phonotnics Program, the National Sciences and Engineering Reasearch Council of Canada via Discovery Grant RGPIN-2023-05818, and the Canada First Research Excellence Fund through its collaboration with the Institut de Valorisation des Données (IVADO).

\end{document}